\documentclass[12pt, legno]{article}
\usepackage[dvips]{graphicx}
\usepackage{pictexwd,dcpic,amssymb,amsmath,amsthm}
\usepackage[mathscr]{eucal}
\begin{document}
\setlength{\hoffset}{-.5in}

\parskip=6pt

\def\tilde{\widetilde}
\def\la{\langle}
\def\ra{\rangle}
\def\fraka{\frak a}
\def\frakz{\frak z}
\def\frakn{\frak n}
\def\var{\varepsilon}
\def\bG{\mathbb G}
\def\bM{\bold M}
\def\bN{\bold N}
\def\bP{\bold P}
\def\bO{\bold O}
\def\bB{\bold B}
\def\bT{\bold T}
\def\bU{\bold U}
\def\bF{\bold F}
\def\bH{\bold H}
\def\cO{\mathcal O}
\def\cP{\mathcal P}
\def\cH{\mathcal H}
\def\calZ{\mathcal Z}
\def\calY{\mathcal Y}
\def\bR{\Bbb R}
\def\bC{\Bbb C}
\def\bQ{\Bbb Q}
\def\ds{\displaystyle}
\def\bA{\Bbb A}
\def\bAA_F{\Bbb A_F}
\def\bI{\Bbb I}
\def\bZ{\Bbb Z}
\def\Gal{\text{Gal}}
\def\diag{\text{diag}}
\def\sym{\text{Sym}}
\def\Hom{\text{Hom}}
\def\Ind{\text{Ind}}
\def\Spin{\text{Spin}}

%${}^L\!H_v$

\title{\bf Arthur Packets and the Ramanujan Conjecture}

\author{By:\ Freydoon Shahidi\thanks{Partially supported by NSF Grant DMS-0700280}\\
To the Memory of Masayoshi Nagata}

\maketitle

\bigskip

\begin{abstract}
The purpose of this paper is to show that under a part of generalized Arthur's $A$--packet conjecture, locally generic cuspidal automorphic representations of a quasisplit group over a number field are of Ramanujan type, i.e., are tempered at almost all places.
The $A$--packet conjecture allows one to reduce the problem to a special case of a general local question about the components of the corresponding Langlands $L$--packet which is then answered here in its generality.
\end{abstract}

\begin{center}\bf 1.\ \ Introduction\end{center}

Among cuspidal representations of a quasisplit reductive group $G$ over a number field are those whose Whittaker Fourier coefficients are non--zero.
Such representations are usually called globally generic (Definition 2.9 via equation (2.8)).
For $GL_n$ every cuspidal automorphic form is such \cite{49}, while for other groups there are examples of cuspidal representations which do not have such non--vanishing coefficients such as Siegel modular forms.

Beside their applications in the theory of automorphic forms and $L$--functions \cite{18,45,48}, generic representations play an extra role in the context of $L$--packets.
These are disjoint subsets of irreducible admissible representations of either $G(k)$ or $G(\bA_k)$, according as $k$ is local or global, respectively.
They are crucial to the analysis of the trace formula \cite{6,7,32,33,36}.
Moreover, members of an $L$--packet are all either tempered or none are \cite{6,32,36}.
The crucial global issue is which members of an $L$--packet are automorphic, i.e., appear in $L^2(G(k)\backslash G(\bA_k))$.

It has been conjectured that, when $k$ is local, every tempered $L$--packet contains a unique generic representation with respect to a fixed generic character of $U(k)$, where $U$ is the unipotent radical of a Borel subgroup of $G$ over $k$ \cite{48}.
Many cases of this conjecture are now proved \cite{15,17,24,30,48,50,51}.
In particular, for every fixed generic character of $U(k)$, the corresponding generic tempered representations are  in one--one correspondence with tempered $L$--packets.
More precisely, they parametrize local tempered $L$--packets and can be used as base points \cite{33}.

On the other hand, starting with the complementary series for $SL_2(k)$, one sees that there are many irreducible generic unitary, but non--tempered representations for $G(k)$, where $k$ is local.
The purpose of this paper is to show that under the validity of a part of Arthur's $A$--packet conjecture (Conjecture 6.1 here) on automorphic representations \cite{2,3,4}, as generalized by Clozel (Conjecture 2A of \cite{m}), this will never happen for {\it automorphic} cuspidal representations of $G(\bA_k)$, where $k$ is a number field.
In other words:\ {\it If Conjecture 6.1 is valid, then locally generic cuspidal automorphic representations of $G(\bA_k)$ are always tempered (Theorem 6.2 and Corollary 6.5)}.
In particular, this shows that the locally generic cuspidal representations of $G(\bA_k)$ obtained by putting together the locally generic ones in each local tempered $L$--packet, if automorphic, exhaust all the automorphic generic representations.
In conclusion, generic representations also parametrize global tempered $L$--packets.

What actually follows from Arthur's conjecture is that almost all the components of a locally generic automorphic form are tempered.
One expects that this extends to all places as we state as Conjecture 6.4 and give reasons for its validity.

Appealing to Arthur's upcoming book \cite{7} (see also \cite{12}) and the automorphic descent of Ginzburg-Rallis-Soudry \cite{18,x}, in (6.8) we sketch an argument for all these conjectures for quasisplit classical groups, under the assumption of validity of the Ramanujan conjecture for $GL_N$ which requires that the cuspidal representations on $GL_N(\bA_k)$ be all tempered.
At present, all that we know for a general $n$, is the estimates in \cite{d}.
Incidentally, our Theorem 6.2 implies that if Conjecture 6.1 is true, then so is the Ramanujan conjecture for $GL_N$ since cusp forms on $GL_N(\bA_k)$ are all generic \cite{49}.
This shows how deep Arthur's $A$--packet conjecture and its generalization are.

There are several consequences of such a result.
First that the general belief \cite{a,43,47} that globally generic cuspidal representations, which are automatically locally generic, are tempered, is in fact correct, if Arthur's conjecture is.
There was no conceptual reason before for why this should be true.
In fact, it was believed that being globally generic, i.e., having a non--vanishing Whittaker Fourier coefficient, is important for the form being tempered, where as our result shows it is being locally generic which matters and is enough.

The second is that there should be no difference between the locally and globally generic representations (Conjecture 2.10).
In fact, all the conjectured restrictions seem to be of relevance only if the representation is not tempered (Conjecture 2.6 of \cite{g}).
(See the discussion at the end of paragraph 6.8 here for classical groups and Conjecture 24.2 of \cite{16}.)

The trace formula is insensitive to detecting globally generic representations and in practice one has to usually use a Poincar\'e series to construct them \cite{48}.
On the other hand, in view of the connection of locally generic cuspidal automorphic representations with tempered $L$--packets and our Theorem 6.2, they may be more amenable to detection by Arthur's trace formula.
In particular, we hope to rule out the existence of locally generic representations which are not generic with respect to local components of any generic character of $U(k)\backslash U(\bA_k)$ by means of multiplicity formulas in \cite{32,c,36}; see also \cite{3,4}.
We plan to take up such questions in a future paper.

The proof relies on an analysis of the Langlands $L$--packet for the parameter $\phi_{\psi_v}$ attached to the Arthur parameter $\psi_v$ for an unramified component of the automorphic representation $\pi=\otimes_v\pi_v$ of $G(\bA_k)$.
We refer to Section 3 for the definitions and basic properties of these parameters.
Using representation theory developed within our method we then show that the parameter $\phi_{\psi_v}$ is tempered if its packet has a generic member (cf.~\cite{10,48} and particularly Proposition 5.4 of \cite{10}).

Our arguments are quite general and are applicable to any $\phi_\psi$ whenever the local Langlands conjecture is valid for proper Levi subgroups of $G$ to the effect of equality of Artin $L$--functions with the corresponding ones defined through the Langlands--Shahidi method for the $k$--points of these Levi subgroups, where $k$ is any local field (Theorem 5.1 and Proposition 4.14).
More precisely, let $G$ be a quasisplit connected reductive group over a local field $k$ and let $P=MN$ be a proper parabolic subgroup of $G$ defined over $k$.
Fix a quasi--tempered irreducible generic representation $\sigma$ of $M(k)$.
Let $r$ be the adjoint action of ${}^LM$, the $L$--group of $M$, on $^L\frak n$, the Lie algebra of the $L$--group of $N$.
Let $L(s,\sigma,r)$ be the $L$--function attached to $\sigma$ and $r$ through the Langlands--Shahidi method defined in \cite{48}.
Our working assumption is then through the local Langlands correspondence (LLC) to the effect that
$$
L(s,r\cdot\phi)=L(s,\sigma,\tilde r),
$$
where the $L$--function on the left is the Artin $L$--function attached to the representation $r\cdot\phi$, if $\phi$ is the homomorphism of the Weil--Deligne group into ${}^LM$ parametrizing $\sigma$.
There are a number of instances that this has already been established \cite{22,25,28,29,43}.
Under this assumption our proof then uses
Proposition 5.4 of \cite{10}, together with the so called ``standard modules conjecture'' and the ``tempered $L$--functions conjecture'' which are now fully proved in \cite{21} (and \cite{20}).
We refer to Section 4 for a discussion of their statements and history.

When $k=\bR$ (or $\bC$), the equality of these factors through the local Langlands correspondence \cite{37} was established in \cite{46}.
Consequently, Corollary 5.3 of our Theorem 5.1 implies that every $\phi_\psi$ whose packet contains a generic member is in fact tempered, whenever $k=\bR$ (or $\bC$).
As it was pointed out to us by Vogan, this may also be proved by other methods in the case of real groups.

In certain cases of classical groups over $p$--adic fields, Theorem 5.1 was proved in \cite{9,39} using a subtle application of classification theorems.
Since the local Langlands conjecture for generic representations to the effect of equality of these $L$--functions is now established (also used in \cite{9,39}) in a number of cases of classical groups \cite{19,22,24}, our theorem also gives a new proof of the aforementioned results in \cite{9,39}.

This is an example of how the parametrization established by the local Langlands conjecture allows us to prove results in representation theory of reductive groups over local fields by means of properties of corresponding Artin $L$--functions.
As our Proposition 4.12 shows, the necessary harmonic analytic data is, in fact, encapsulated in the poles of these $L$--functions.

The paper is organized as follows.
Basic definitions and conjectures are discussed in Section 2.
Arthur parameters $\psi$, their attached Langlands parameters $\phi_\psi$ and their properties are explained in Section 3.
If $\psi=(\phi,\rho)$, where $\rho$ is the $SL_2(\bC)$--component of $\psi$, then Proposition 3.30 gives a characterization of the image of $\rho$ if $\rho\neq 1$, i.e., if $\psi$ is not tempered, which is essential to the proof of Theorem 4.1 which covers the unramified cases.
The general case is stated and proved as Proposition 3.31.

Section 4 discusses the main representation theoretic tools that we need in order to prove our main local Theorem 5.1 from which Theorem 4.1 (Corollary 5.4) follows as a special case.
These results and techniques are mainly extracted from our method since the representations involved are generic.
The main tool here is Proposition 5.4 of \cite{10} and the two conjectures now completely proved in \cite{21}.
Here they are stated as Theorem 4.2 and statement (4.3).

Theorem 5.1 is then proved by putting together Proposition 3.31 and the material in Section 4, particularly Proposition 4.14.
The global consequences and our main global result, Theorem 6.2, is then proved in Section 6.
We finish the section by giving more evidence for each of conjectures involved.

{\bf Acknowledgements}:\ I would like to thank James Arthur for a number of helpful correspondences and for patiently answering my questions.
Thanks are also due to James Cogdell and Dihua Jiang as well as Dick Gross, Dipendra Prasad and Dinakar
Ramakrishnan for useful communications on Conjecture 2.10.
I would also like to thank Sandeep Varma for his careful reading of the paper and many comments, particularly on Conjecture 6.4.
My interest in studying Arthur packets was renewed after I was asked to give some lectures at the Atlas Project Summer Graduate School at University of Utah in July of 2009.
I would like to thank the organizing commitee, especially Jeff Adams and David Vogan for their invitation.
Thanks are also due to Vogan for useful discussions and communications.
Finally, it is my pleasure to dedicate this paper to the Memorial Volume in honor of Professor Masayoshi Nagata for which I would like to thank the Editorial Committee for their invitation.

Thanks are due to the referee for a careful reading of the paper and many comments on both its substance and exposition.

\vfill\eject
\begin{center}\bf 2.\ \ Preliminaries\end{center}

Let $k$ be a number field and denote by $\bA_k$ its ring of adeles.
For each place $v$ of $k$, we let $k_v$ be the completion of $k$ at $v$.
Let $\cO_v$ and $\cP_v$ be the ring of integers and its maximal ideal, respectively.
We let $\varpi_v$ be a generator of $\cP_v$ and normalize an absolute value so that $|\varpi_v|=q_v^{-1}$, where $q_v$ is the cardinality of the field $\cO_v/\cP_v$.

Let $G$ be a quasisplit connected reductive group over $k$.
We fix a Borel subgroup $B$ over $k$ and write $B=TU$, where $T$ is a maximal $k$--torus of $G$ isomorphic to the quotient $B/U$, where $U$ is the unipotent radical of $B$.

We will then have $G(k_v)$ and $G(\bA_k)$ for the $k_v$ and $\bA_k$--points of $G$ as well as for any of subgroups of $G$ defined over $k$.

We remark that for almost all $v,\ G$ is defined over $\cO_v$ and splits over an unramified extension of $k_v$.
Then $G(\cO_v)$ is a hyperspecial maximal compact subgroup of $G(k_v)$.

The choice of the Borel subgroup defines a set of positive roots of $G$, i.e., roots of $T$ on Lie$(U)$.
We let $\tilde\Delta=\tilde\Delta(T,G)$ be the set of simple (non--restricted) roots among them.
Let $\{X_\alpha|\alpha\in\tilde\Delta\}$ be a choice of root vectors in Lie$(U)$.
This means that there exists a map
\[
\phi\colon U\longrightarrow\prod\bG_a,\tag{2.1}
\]
where the product runs over all the roots in $\tilde\Delta$, sending $\exp(x_\alpha X_\alpha)$ to $x_\alpha,x_\alpha\in\overline k$, whose kernel contains the derived group of $U$.
Composing $\phi$ with the map
\[
\Sigma\colon\prod_{\alpha\in\tilde\Delta}\bG_a\longrightarrow\bG_a,\tag{2.2}
\]
defined by
\[
\Sigma((x_\alpha)_\alpha)\mapsto\sum_{\alpha\in\tilde\Delta} x_\alpha,\tag{2.3}
\]
we then get a map $\Sigma\cdot\phi$ from $U$ to $\bG_a$.
If the Galois group $\Gamma=\Gal(\overline k/k)$ fixes the splitting $\{X_\alpha\}_\alpha$ as a set, then $\Sigma\cdot\phi$ is defined over $k$ and the splitting is said to be defined over $k$.
For a quasisplit group such splittings exist which we will fix one from now on.
This definition is valid for $k$ as well as each $k_v$ and in fact a splitting over $k$ will be also one over $k_v$ for each $v$.

According as $k$ is local or global, we fix a non--trivial character $\psi_v$ or $\psi$ of $k_v$ or $k\backslash\bA_k$, respectively.
We then define a generic character $\chi_v$ or $\chi$ of $U(k_v)$ or $U(k)\backslash U(\bA_k)$ by
\[
\chi=\psi\cdot\Sigma\cdot\phi,\tag{2.4}
\]
respectively.
When $k$ is global, then $\psi=\otimes_v\psi_v$ and thus $\chi=\otimes_v\chi_v$, where $\chi_v|U(\cO_v)=1$ for almost all places $v$.

Now, let $\pi=\otimes_v\pi_v$ be a cuspidal automorphic representation of $G(\bA_k)$.
We shall assume $\pi$ is unitary.
This simply means that its central character $\omega_\pi$ is unitary.

A representation $\sigma$ of $G(k_v)$ on a complex vector space $\cH(\sigma)$ is called $\chi_v$--generic, if there exists a functional $\lambda_v$ on the continuous dual $\cH(\sigma)'$ of $\cH(\sigma)$ such that
\[
\lambda_v(\sigma(u) w)=\chi_v (u) \lambda_v (w)\tag{2.5}
\]
for every $w\in\cH(\sigma)$.
When $k_v$ is archimedean, one requires the continuity to be with respect to the semi--norm topology on the space of differentiable vectors $\cH(\sigma)_\infty$, the span of $\sigma(f)w$, $f\in C_c^\infty (G(k_v))$, $w\in\cH(\sigma)$.
For this, one requires the topology on $\cH(\sigma)$ be defined by either Banach or more generally Frech\'et space norms or semi--norms, respectively.
In particular, the Hilbert norm may be used when $\sigma$ is unitary.
It is well--known (cf.~\cite{49}) that if $\sigma$ is irreducible, then the space of such functionals is at most one--dimensional.

If $\sigma$ is $\chi_v$--generic, then one can fix a $\chi_v$--Whittaker functional $\lambda_v$ on $\cH(\sigma)$ and for each vector $x\in\cH(\sigma)_\infty$, define a Whittaker function $W_x(g)$ on $G(k_v)$ by
\[
W_x(g)=\lambda_v (\sigma(g)x).\tag{2.6}
\]
Up to a complex multiple, the model is unique.

\medskip\noindent
{\bf Definition 2.7}.
{\it A cuspidal representation $\pi=\otimes_v\pi_v$ is called locally generic if each $\pi_v$ is generic with respect to a generic character $\chi_v$ of $U(k_v)$.
Note that we are not requiring $\chi_v$ to be a local component of a global character $\chi$ of $U(k)\backslash U(\bA_k)$}.
\medskip

There is also a notion of a globally generic cusp form.
Let $\chi$ be a generic character of $U(k)\backslash U(\bA_k)$.
Assume $\pi$ is a cuspidal representation of $G(\bA_k)$ and let $\varphi$ be a cusp form in the space of $\pi$.
Let
\[
W_\varphi(g)=\int_{U(k)\backslash U(\bA_k) }\varphi(ug)\overline{\chi(u)}\ du,\tag{2.8}
\]
which converges for all $\varphi$.
We then define

\medskip\noindent
{\bf Definition 2.9}.
{\it A cuspidal representation $\pi=\otimes_v\pi_v$ is called globally generic with respect to $\chi=\otimes_v\chi_v$, if $W_\varphi (e)\neq 0$ for some $\varphi\in\cH(\pi)$.
Note that $\pi$ is then automatically locally generic}.
\medskip

Being globally generic seems to be a much stronger condition.
One goal of this paper is to provide enough evidence for the following conjecture.

\medskip\noindent
{\bf Conjecture 2.10}.{\it Assume $\pi=\otimes_v\pi_v$ is locally generic with respect to local components of a generic character $\chi=\otimes_v\chi_v$ of $U(k)\backslash U(\bA_k)$.
Then $\pi$ is globally generic.
More precisely, the isotypic constituent of the space of cusp forms containing $\pi$ is globally generic, i.e., (2.8) is non--vanishing on it.}

\medskip\noindent
{\it Remark 2.11}.
The conjecture is a well--known theorem
for $G=GL_n$, since all the cusp forms on $GL_n(\bA_k)$ are globally generic \cite{49}, and it seems to agree with Conjecture 24.2 of \cite{16} as well as the discussion in Paragraph (6.8) here.
It is true even when multiplicity one fails (e.g.~\cite{e,f}), which could happen as predicted by the trace formula (\cite{3,33}).
(See Remark 2.13.)
On the other hand there are many examples of non--generic cuspidal representations for other groups 
\cite{35,44}.
Among them are the so called CAP representations:\ {\it Cuspidal Representations Associated to Parabolics} as coined by Piatetski--Shapiro.
They include the examples of Saito--Kurokawa \cite{b} and Howe--Piatetski--Shapiro \cite{a}.
What is special about CAP representations is that none are tempered.
As explained in the introduction, one may look at the problem in terms of $L$--packets and here is where one can expect that generic representations completely parametrize global tempered $L$--packets.
What is surprising is that the converse seems to be also true.
In fact, in this paper, using a part of a conjecture of Arthur on global $A$--packets and a natural rigidity conjecture, we prove (Theorem 6.2 and Corollary 6.5):

(2.12)\ {\it Assume the validity of Conjectures 6.1 and 6.4.
Then locally generic cuspidal automorphic representations for $G(\bA_k)$ are tempered}.

This is purely global.
There are irreducible generic unitary spherical representations of even $GL_2(k_v)$ which are not tempered; the complementary series.

Statement (2.12) guarantees that there are no global obstructions for the equivalence of locally and globally generic cuspidal representations, at least up to isomorphisms, since these obstructions are automatically satisfied for tempered representations.
These obstructions are usually stated in terms of adjoint $L$--functions for the cuspidal representation.
Note that isomorphic irreducible admissible representations of $G(\bA_k)$ have same $L$--functions.

\medskip
{\it Remark 2.13}.
Assume multiplicity one fails for a cuspidal generic representation $\pi$ of $G(\bA_k)$.
Let $A_1$ and $A_2$ be two linearly independent embeddings of $\pi$.
Let $\lambda$ be the global Whittaker functional on the space of cusp forms, i.e.,
$$
\lambda(\varphi)=\int_{U(k)\backslash U(\bA_k) }\varphi (ug)\overline\chi (u) du
$$
as in (2.8).
Choose $c\neq 0$ such that
$$
\lambda\cdot A_1=c\lambda\cdot A_2.
$$
Then $\lambda\cdot (A_1-cA_2)=0$.
Thus $A_1-cA_2$ provides a (non--zero) embedding of $\pi$ which is not globally generic.
(We were reminded of this construction by D.~Prasad.)
The content of Conjecture 2.10 is that $\pi$ still appears as a globally generic representations.
In particular, any locally generic cuspidal representation should also embed as a globally generic one, and therefore the theory of $L$--functions developed through different methods for globally generic representations is valid whether $\pi$ embeds as a globally or locally generic representation.

\begin{center}3. \ \ {\bf Arthur Parameters}\end{center}

We continue to assume that $G$ is quasisplit as this will greatly reduce the notation and is sufficient for the purposes of this paper.

We first assume $k$ is local.
Let $L_k$ be either $W'_k$, the Weil--Deligne group of $k$, if $k$ is non--archimedean, i.e., $W'_k=W_k\times SL_2(\bC)$ or the Weil group, otherwise.
Let $\Phi(G)=\Phi(G/k)$ be the set of Langlands parameters, i.e., equivalence classes of homomorphisms 
$$
\phi\colon L_k\longrightarrow {}^L\!G=\widehat G\rtimes L_k
$$
under conjugation by elements in $\widehat G$, satisfying the following conditions:

(3.1) $\phi$ factors through $L_{K/k}$ for some finite extension $K/k$, where $L_{K/k}$ is the corresponding Weil--Deligne or Weil group defined by $K/k$, i.e., $L_{K/k}=W_{K/k}\times SL_2(\bC)=W'_{K/k}$ or $L_{K/k}=W_{K/k}$ according as $k$ is non--archimedean or not; 

(3.2) $\phi$ is continuous on $W_{K/k}$ and complex analytic on $SL_2(\bC)$;

(3.3) the diagram
\[\begindc{0}[50]
\obj(1,0){$L_{K/k}$}
\obj(3,1){$^LG=\widehat G\rtimes L_k$}
\obj(0,1){$L_k$}
\mor(0,1)(2,1){$\phi$}
\mor(0,1)(1,0){}
\mor(2,1)(1,0){$\text{Proj}_2$}
\enddc\]
commutes;

(3.4) the image of $W_{K/k}$ consists of only semisimple elements.

We note that (3.1) and (3.2) mainly address the continuity of $\phi$.

We finally let $\Phi_{\text{temp}}(G)=\Phi_{\text{temp}}(G/k)$ denote those $\phi\in\Phi(G)$ whose image in $\widehat G$ is bounded.

We now define Arthur parameters.
Let $\Psi(G)=\Psi(G/k)$ be the set of $\widehat G$--orbits of maps
\[
\psi\colon L_k\times SL_2 (\bC)\longrightarrow {}^LG=\widehat G\rtimes L_k\tag{3.5}
\]
such that the projection of $\psi (L_k)$ onto $\widehat G$ is bounded.
Moreover, we assume $\phi=\psi|L_k\in\Phi_{\text{temp}}(G)$, i.e., that in addition it satisfies (3.1)--(3.4).

Finally for each such $\psi\in\Psi(G)$, define a Langlands parameter $\phi_\psi\in\Phi(G)$ by
\[
\phi_\psi(w)=\psi(w,\begin{pmatrix} |w|^{1/2}&0\\ 0&|w|^{-1/2}\end{pmatrix}).\tag{3.6}
\]

We now recall an argument from \cite{2} which implies that:

(3.7) {\it The map
$$
\psi\mapsto\phi_\psi
$$
is an injection from $\Psi(G)$ into $\Phi(G)$.}

By Jacobson--Morozov (and Kostant) there is a one--one correspondence between unipotent conjugacy classes in a complex semisimple (or reductive) group and conjugacy classes of maps from $SL_2(\bC)$ to the group (cf.~\cite{14}).

Each such conjugacy class of unipotent elements then determines a ``distinguished'' semisimple conjugacy class (cf.~\cite{14}) as the orbit of the semisimple member $H$ of the $SL_2(\bC)$--triple which gives the map.

One then attaches a ``weighted'' Dynkin diagram to each such orbit by numbering the simple roots $\alpha$ by $\alpha(H)$, where $H$ is dominant with respect to the set of simple roots.
It can be shown that $\alpha(H)\in \{0,1,2\}$ and that weighted Dynkin diagrams are in 1--1 correspondence with unipotent conjugacy classes or conjugacy classes of maps from $SL_2(\bC)$ into the group (cf.~\cite{14}).

Now each $\psi\in\Psi(G)$ can be decomposed as (cf.~\cite{2})
\[
\psi(w,g)=\phi(w)\rho(g),\tag{3.8}
\]
where $\phi=\psi|L_k\in\Phi_{\text{temp}}(G)$ and if
\[
C_\phi=\text{Cent}(\text{Im}(\phi),\widehat G),\tag{3.9}
\]
then
\[
\rho\colon SL_2 (\bC)\longrightarrow C_\phi.\tag{3.10}
\]

By the earlier discussion, weighted Dynkin diagrams, which are defined by the distinguished orbits, are in one--one correspondence with unipotent conjugacy classes in $\widehat G$, or with orbits of maps from $SL_2(\bC)$ to $\widehat G$.
In particular, for $\rho$ as in (3.10), the restriction of $\rho$ to
$$
\left\{\begin{pmatrix} t&0\\ 0&t^{-1}\end{pmatrix}\bigg| t\in\bC^*\right\}
$$
determines the conjugacy class of $\rho$.
Thus $\phi_\psi$ determines $\psi$, and thus the injection in (3.7).

Next, we recall that for each such $\psi$ over a local field $k$, Arthur conjectured the existence of a finite set $\Pi(\psi)$ of irreducible admissible representations of $G(k)$, satisfying a number of properties (cf.~\cite{2,3,4}).
In particular, he demanded that the $L$--packet $\Pi(\phi_\psi)$ defined by the Langlands parameter $\phi_\psi$ be included in $\Pi(\psi)$, i.e.,
\[
\Pi(\phi_\psi)\subset\Pi(\psi).\tag{3.11}
\]
Note that while the members of $\Pi(\psi)$ are rather mysterious and are there to supplement those in $\Pi(\phi_\psi)$ to produce suitable stable distributions (cf.~\cite{2,3,4}), those in $\Pi(\phi_\psi)$ are readily available through Langlands classification \cite{37} so long as the local Langlands conjecture (LLC) is known for the ``defining'' Levi subgroup $M$ for $\Pi(\phi_\psi)$ as we now review.
We recall (cf.~\cite{2,3}) that a member $\phi$ of $\Phi(G)$ always determines a commuting pair $\phi_0\in\Phi_{\text{temp}}(G)$ and a $\phi_+\in\Phi(G)$ such that
\[
\phi(w)=\phi_0(w)\phi_+(w)\qquad (w\in L_k),\tag{3.12}
\]
and so that $\phi\in\Phi_{\text{temp}}(G)$ whenever $\phi_+=1$.
The centralizer of the image of $\phi_+$ in ${}^LG$ is a Levi subgroup ${}^LM$.
The subgroup $^LM$ is what we have chosen to call the ``defining'' Levi subgroup of the packet $\Pi(\phi)$.
We note that it is possible to have $^LM={}^LG$.

We point out that when $k=\bR$, candidates for Arthur packets have been proposed for any connected reductive real group in \cite{1}.
All the conjectured properties of these packets are also verified in \cite{1}, except for the most difficult:

(3.13) {\it Every member of each Arthur packet $\Pi(\psi)$, $\psi\in\Psi(G/\bR)$, is unitary}.

For $p$--adic groups, except for the upcoming work of Arthur on classical groups, these packets are known only sporadically.
Still one can make some assertions which have global consequence.

Let $I'_k\subset L_k$ be $I'_k=I_k\times SL_2(\bC)$, where $I_k$ is the inertia subgroup of $W_k$.
Fix $\psi\in\Psi(G)$.
Decompose $\psi$ as in (3.8):\ $\psi=(\phi,\rho)$.
Assume $\phi|I'_k\equiv 1$.
Then $\phi_\psi|I'_k\equiv 1$ and conversely.
Consequently $\Pi(\phi_\psi)$ will consist of unramified representations of $G(k)$ each for a hyperspecial maximal compact subgroup of $G$.
We recall that they will constitute all the unramified constituents of a principal series whose inducing character is parametrized by $\phi_\psi$ which up to conjugation now factors through ${}^LT$:
\[
\phi_\psi\colon L_k\longrightarrow {}^LT.\tag{3.14}
\]

If $G$ is defined over $\cO=\cO_k$, the ring of integers of $k$, then $\Pi(\phi_\psi)$ will have a unique unramified representation with respect to $G(\cO)$.

Now assume $k$ is a number field.
We will then use $L_k$ to denote the hypothetical global Langlands group \cite{5,32}.
Then there should be maps
\[
L_{k_v}\longrightarrow L_k.\tag{3.15}
\]
We will again let $\psi(G)=\psi(G/k)$ be the set of equivalence classes of maps
\[
\psi\colon L_k\times SL_2 (\bC)\longrightarrow {}^LG\tag{3.16}
\]
for which the image of $L_k$ in $\widehat G$ is bounded.
This time equivalence between two maps $\psi_i$, $i=1,2$, means that there exists an element $g\in\widehat G$ such that
\[
g^{-1}\psi_1(w,x)g=\psi_2(w,x)z_w,\tag{3.17}
\]
where $z_w$ is a 1--cocycle of $L_k$ in the center $Z(\widehat G)$ of $\widehat G$ whose class in $H^1(L_k,Z(\widehat G))$ is locally trivial.
This then agrees with local equivalence.

We can then consider the map
$$
\psi\mapsto\phi_\psi
$$
from $\Psi(G/k)$ into $\Phi(G/k)$, the set of global Langlands parameters (cf.~\cite{4}).

The maps from each $L_{k_v}$ to $L_k$ then allow us to define $\psi_v\in\Psi(G/k_v)$ and $\phi_{\psi_v}\in\Phi(G/k_v)$.

Given $\psi\in\Psi (G/k)$ we define the global Arthur packet to be
\[
\Pi(\psi)=\{\pi=\otimes_v\pi_v|\pi_v\in\Pi(\psi_v)\},\tag{3.18}
\]
where for almost all $v$, $\pi_v=\pi_v^0$, the unique $G(\cO_v)$--spherical representation in $\Pi(\phi_{\psi_v})$, the $L$--packet attached to $\phi_{\psi_v}$.

Arthur's conjecture then states that every automorphic representation must belong to $\Pi(\psi)$ for some $\psi\in\Psi(G/k)$.

For simplicity of notation we shall now assume $k$ is local and drop the subscript $v$.

Given the parameter $\phi_\psi$ attached to a $\psi\in\Psi(G/k)$, the pair $\phi_0$ and $\phi_+$ in (3.12) are $\phi$ and 
\[
\phi_+\colon w\mapsto\rho (\begin{pmatrix} |w|^{1/2}&0\\ 0&|w|^{-1/2}\end{pmatrix} ),\tag{3.19}
\]
respectively, with $\phi$ and $\rho$ as in (3.8).
Both $\phi$ and $\phi_+$ have their images in ${}^LT$ upon conjugation.

If $\rho\neq 1$, then there exists a positive root $\widehat\alpha$ of $\widehat T$ in Lie$(\widehat U)$ such that
for each $t>0$
\[
\log_t (\hat\alpha (\rho\begin{pmatrix} t^{1/2}&0\\ 0&t^{-1/2}\end{pmatrix}))=1,\tag{3.20}
\]
\[
while\; the\; adjoint\; action\; of \;\phi(L_k)\; on\; X_{\widehat\alpha}\; is\; trivial,\tag{3.21}
\]
since the nilpositive element $X$ of the $sl_2(\Bbb C)$--triple attached to $\rho$ lies
in the centralizer $C_\phi$ of Im$(\phi)$.
We will then say $X_{\widehat\alpha}$ contributes to the nilpositive element $X$.
This simply means that if $X$ is written as a linear combination of $X_{\widehat\beta},\widehat\beta>0$, which is possible upon conjugation, then the coefficient of $X_{\widehat\alpha}$ is non--zero.

Next, let
\[
\mu\colon T(k)\longrightarrow \bC^*\tag{3.22}
\]
be the character of $T(k)$ attached to $\phi$ by Langlands \cite{34,38}.
Similarly, let
\[
\nu\colon T(k)\longrightarrow \bC^*\tag{3.23}
\]
be the one attached to $\phi_+$.

Let $\varpi$ be a uniformizer in $k$ and denote by $H_\gamma$ the coroot at a root $\gamma$.
The simple roots $\gamma$ of $A_0$, the split component of $T$, in Lie$(U)$ for which
\[
|\nu (H_\gamma (\varpi))|=1\tag{3.24}
\]
generate a Levi subgroup $M$ of $G$, $M\supset T$, whose connected $L$--group $\widehat M={}^L\!M^0$ contains Im$(\phi)$ in $\widehat G$.

Let $P=MN$ be the parabolic subgroup of $G$ with $M$ as a Levi subgroup and $N\subset U$.
Let $U_M=U\cap M$.
Then the representation
\[
\tau=\Ind^{M(k)}_{T(k) U_M(k)} \mu\tag{3.25}
\]
is a tempered representation of $M(k)$, which may not be irreducible.
Write
\[
\tau=\bigoplus^n_{i=1} \tau_i\tag{3.26}
\]
for a decomposition of $\tau$ to its irreducible subrepresentations.
We first observe:

(3.27) The root $\alpha$ for which $X_{\widehat\alpha}$ contributes to the nilpositive element in the
$sl_2(\bC)$--triple attached to $\rho$ may be chosen not to be
in Lie$(U_M)$ since the restriction of $\hat\alpha\cdot\rho$ to diagonal elements in $SL_2(\bC)$ equals $\widehat\alpha\cdot\phi_+$ which is non--trivial.

Next, we note that replacing $\nu,\tau_i$ and $M$ with a $W(G,A_0)$--conjugate, we may assume that
\[
I(\nu,\tau_i)=\Ind^{G(k)}_{P(k)}\ \tau_i\otimes\nu\tag{3.28}
\]
is in the Langlands setting and is thus a standard module.
Assertion (3.27) still remains valid.
The members of $\Pi(\phi_\psi)$ are now the unique Langlands quotients $J(\nu,\tau_i)$ of each $I(\nu,\tau_i)$, $1\leq i\leq n$, i.e.,
\[
\Pi(\phi_\psi)=\{J(\nu,\tau_i)\ |\ 1\leq i\leq n\}.\tag{3.29}
\]
This is part of the local Langlands conjecture (LLC) (cf.~\cite{37}, for example) which is automatic for unramified representations.

We have therefore shown:

\medskip\noindent
{\bf Proposition 3.30}.
{\it Assume $k$ is a non--archimedean local field.
Fix $\psi\in\Psi(G/k)$ and consider the Langlands parameter $\phi_\psi$ attached to $\psi$.
Decompose $\psi=(\phi,\rho)$ as in (3.8).
Assume $\phi$ is trivial on $I'_k$.
Let $\Pi(\phi_\psi)$ be the $L$--packet attached to $\phi_\psi$.
Assume $\rho\neq 1$.
Then
$$
\Pi(\phi_\psi)=\{J(\nu,\tau_i)| 1\leq i\leq n\},
$$
where each $\tau_i$ is an irreducible tempered representation of a proper Levi subgroup $M(k)$ of $G(k)$ for which there exists a root $\alpha$ with
$X_\alpha\!\in \!\text{Lie}(N(k))$, $N\subset U$, where $X_{\widehat\alpha}\!\in\!\text{Lie}(\widehat U)$ contributes to
the nilpositive member of the $sl_2(\bC)$--triple defined by $\rho$.
Moreover
$$
\log_q ( \widehat\alpha (\rho \begin{pmatrix} q^{1/2}&0\\ 0&q^{-1/2}\end{pmatrix} ))=1,
$$
while the adjoint action of $\phi(L_k)$ on $X_{\widehat\alpha}$ is trivial,
i.e., (3.20) and (3.21) are valid with $t=q$, where $q$ is the number of elements in the residue field}.

\medskip\noindent
{\it Remark 3.31}.
We point out that although the representations in $\phi_\psi$ are all spherical, their Langlands data is not necessarily induced from a Borel subgroup.
The fact that $X_{\widehat\alpha}\in\text{ Lie}(\widehat N)$ and not only Lie$(\widehat U)$, is crucial for the proof of our main result.

Let us now consider the general case and thus remove the condition $\phi|I'_k=1$.
Write $\phi_\psi=\phi_0\phi_+$ as before.
Again $\phi_0=\phi$ and $\phi_+$ is as in (3.19).
Let ${}^L\!M$ be the centralizer of Im$(\phi_+)$ in ${}^LG$, a Levi subgroup of ${}^LG$.
Let $M$ be the corresponding Levi subgroup of $G$ which we may assume, upon conjugation, to contain $T$.
The parameter $\phi_+$ will be a map into ${}^LT$ and in fact ${}^L\!A$, where $A$ is the split component of $M$.
It defines a positive quasi--character $\nu$ of $M(k)$.
Moreover, $\phi\in\Phi_{\text{temp}}(M/k)$.
Let $\Pi_\phi^M =\Pi^M(\phi)$ be the $L$--packet of $M(k)$ attached to the tempered parameter $\phi$ which we assume exists by assuming LLC for $M$.
This is always the case when $k=\bR$ (or $\bC$) \cite{51} and is also true in certain special cases when $k$ is $p$--adic.
We will discuss this assumption further in later sections.

The $L$--packet $\Pi(\phi_\psi)$ of $G(k)$ attached to $\phi_\psi$ then consists of Langlands quotients $J(\nu,\tau)$ as $\tau$ ranges in $\Pi_\phi^M$.
Again (3.20) and (3.21) are valid.
We therefore have:

\medskip\noindent
{\bf Proposition 3.31}.
{\it Assume $\psi\in\Psi (G/k)$ is arbitrary and choose $\phi,\rho,\phi_+,M$ and $\nu$ as above.
Assume $\Pi^M_\phi$ is defined.
Then by LLC
$$
\Pi(\phi_\psi)=\{J(\nu,\tau) |\tau\in\Pi_\phi^M \},
$$
where
$$
I(\nu,\tau)=\text{Ind }_{M(k) N(k) }^{G(k) }\quad  (\tau\otimes\nu)\otimes\bold 1
$$
are the corresponding standard modules with $N\subset U$.
Assume $\rho\neq 1$.
Then there exists a root $\alpha$ with
$X_\alpha\!\in\text{Lie}(N(k))$, where $X_{\widehat\alpha}\in$ Lie$(\widehat U)$ contributes to the nilpositive member of the $sl_2(\bC)$--triple defined by $\rho$.
Moreover (3.20) and (3.21) are valid}.

\vfill\eject
\begin{center}\bf 4. Some Representation Theory\end{center}

In this section we review the main tools from representation theory that we need in order to prove the main step in the proof of (2.12), namely:

\medskip\noindent
{\bf Theorem 4.1}.
{\it If a member of the $L$--packet defined by $\phi_\psi$ is unramified and generic, then it is tempered}.

The unramified condition in Theorem 4.1 can be removed whenever the local Langlands holds for the Levi subgroup defined by the tempered parameter $\phi$ as we explain later.

To proceed, we assume $k$ is a local field, archimedean or non--archimedean.
Our group $G$ continues to be quasisplit over $k$ with a Borel subgroup $B=TU$ defined over $k$.
We fix a generic character $\chi$ of $U(k)$ by means of our splitting as in Section 2.

Let $I(\nu,\sigma)$ be a standard module for $G(k)$ defined by a parabolic subgroup $P=MN$, where $N\subset U$ and $T\subset M$.
Then $\nu$ is in the Langlands setting, i.e., $\nu$ is in the positive Weyl chamber and $\sigma$ is an irreducible tempered representation of $M(k)$.
We let $J(\nu,\sigma)$ be the corresponding (unique) Langlands quotient.
It follows from Rodier's theorem that $J(\nu,\sigma)$ is $\chi$--generic only if $\sigma$ is $\chi_M$--generic, where $\chi_M=\chi|U_M(k)$ with $U_M=U\cap M$.

There is more to be said.
The following theorem is now completely proved in the case of a $p$--adic field $k$ in a recent manuscript of Heiermann and Opdam \cite{21}.

\medskip\noindent
{\bf Theorem 4.2}.
{\it $J(\nu,\sigma)$ is $\chi$--generic if and only if $I(\nu,\sigma)$ is irreducible and $\sigma$ is $\chi_M$--generic}.

This is true for any local field and when $k=\bR$ (or $\bC$) was first proved by Vogan more than 30 years ago in
\cite{51}, using a characterization of generic representations by Kostant \cite{31}.

Before \cite{21}, the $p$--adic case was proved in many instances in a number of papers by other authors \cite{8,10,20,26,27,29,42,48}.
Its proof is reduced to another conjecture, sometimes called the ``tempered $L$--function conjecture'' which was conjectured in \cite{48} and demands that:

(4.3) {\it All the $L$--functions $L(s,\sigma,r_i)$ defined in \cite{48} are holomorphic for Re$(s)~>~0$ whenever $\sigma$ is tempered}.

The progress on this conjecture followed the same path for which one can refer to the references above.

This conjecture, which is now a theorem in \cite{21} as well, is also needed in our argument and we will therefore quickly review both of them.

We start by recalling what a local coefficient is.
For simplicity let us consider only the case of a standard module $I(\nu,\sigma)$, although these objects can be defined very generally.

Let $W(G,A_0)$ be the Weyl group of $A_0$ in $G$.
If $M$ is defined by a subset $\theta$ of simple roots $\Delta$ of $A_0$ in $U$, we let $\tilde w$ be an element in $W(G,A_0)$ such that $\tilde w(\theta)\subset\Delta$.
Although it is not relevant here, we let $w$ be a representative for $\tilde w$ as in \cite{46,48}.

We let $A(\nu,\sigma,w)$ be the standard intertwining operator
\[
A(\nu,\sigma,w)f(g)=\int_{N_{\tilde w}}f(w^{-1} ng) dn\tag{4.4}
\]
from $I(\nu,\sigma)$ to $I(w(\nu),w(\sigma))$, where $f$ is in the space of $I(\nu,\sigma)$, $g\in G(k)$ and $N_{\tilde w}=U\cap w\overline N w^{-1}$ with $\overline N$ the unipotent subgroup opposed to $N$.

Next, let $w_0=w_\ell\cdot w^{-1}_{\ell,\theta}$, where $w_\ell$ and $w_{\ell,\theta}$ are the representatives for the long elements $\tilde w_\ell$ and $\tilde w_{\ell,\theta}$ of $W(G,A_0)$ and $W(M,A_0)$, respectively.

We now assume $\sigma$ is generic with respect to $\chi_M$.
Let $\lambda_M$ be a Whittaker functional for $\sigma$, i.e., so that
\[
\lambda_M(\sigma(u)v)=\chi_M (u) \lambda_M (v),\tag{4.5}
\]
$u\in U_M(k)$, $v\in\cH(\sigma)$.
The functional is also supposed to be appropriately continuous when $k=\bR$ (cf.~\cite{49}).
Let $M'=M_{\tilde w_0(\theta)}$ and let $P'=M'N'$ be the corresponding standard parabolic subgroups of $G$ and thus $N'\subset U$.
The induced Whittaker functional is then defined by
\[
\lambda (f)=\int_{N'(k)} \lambda_M (f(w_0^{-1} n'))\chi^{-1} (n') dn'.\tag{4.6}
\]
We define similarly the induced functional $\lambda'$ for $I(w(\nu),w(\sigma))$.
The local coefficient $C_\chi(\nu,\sigma,w)$ is the complex number defined by
\[
C_\chi(\nu,\sigma,w) \lambda' (A(\nu,\sigma,w)f)=\lambda(f).\tag{4.7}
\]
It is a meromorphic function of $\nu$ which becomes a product of rational functions of one variable each given in the exponentials defined by $\nu$ (cf.~\cite{45,48}), if $k$ is $p$--adic.
When $k=\bR$, $C_\chi(\nu,\sigma,w)$ becomes a ratio of products of $\Gamma$--functions \cite{46}.

Poles and zeros of this function play a crucial role in determining the reducibility of $I(\nu,\sigma)$ for a generic $\sigma$ (cf.~\cite{10,45,48}).

For our purposes and in view of Theorem 4.2, we need to determine the reducibility of the standard module $I(\nu,\sigma)$ and for that we will use Proposition 5.4 of \cite{10} which we now recall as 

\medskip\noindent
{\bf Proposition 4.8 (cf.~\cite{10})}.
{\it Let $k$ be a local field of characteristic zero and let $I(\nu,\sigma)$ be a standard module for $G(k)$, where $\sigma$ is a $\chi_M$--generic representation of $M(k)$.
Then $I(\nu,\sigma)$ is irreducible if and only if $C_\chi (\nu,\sigma,w_0)^{-1}\neq 0$.}

There are a number of cases for which $C_\chi(\nu,\sigma,w_0)$ is computed in terms of the Langlands parameter for $\sigma$.
For example, when $k=\bR$ (or $\bC$) the local coefficients are explicitly computed in terms of Artin $L$--functions in \cite{37,46}.
More precisely, let $\sigma_\nu=\sigma\otimes\exp\langle\nu,H_M(\cdot)\rangle$, where 
$$
H_M\colon M(\bR)\longrightarrow 
\frak a=\Hom(X(M)_{\bR},\bR)$$
is the natural map.
Let
\[
\phi_\nu\colon L_k\longrightarrow {}^L\!M\tag{4.9}
\]
denote the Langlands parameter for $\sigma_\nu$ in $\Phi(G/\bR)$.
Next, let ${}^L{\frak n}$ be the Lie algebra of ${}^L\!N$, the $L$--group of $N$ and denote by $r$ the adjoint action of ${}^L\!M$ on ${}^L\frak n$.
Then
\[
C_\chi(\nu,\sigma,w_0)^{-1} \sim L(0,r\cdot\phi_\nu)/L(1,\tilde r\cdot\phi_\nu),\tag{4.10}
\]
where $\tilde r$ is the contragredient of $r$, which can be realized as the adjoint action of ${}^L\!M$ on 
${}^L{\overline{\frak n}}$, the Lie algebra of ${}^L\!\overline{N}$.
Here $\sim$ signifies the equivalence up to an exponential in $\nu$.

As explained earlier, since $\nu$ is in the positive Weyl chamber and $r\cdot\phi$ is a unitary representation, $\sigma$ being tempered,
\[
L(0,r\cdot \phi_\nu)^{-1}\neq 0\tag{4.11}
\]
by (4.3) (or its Artin $L$--function version \cite{19,22,23,24,25}).
Therefore the zeros of $L(1,\tilde r\cdot\phi_\nu)^{-1}$ will become precisely those of $C_\chi(\nu,\sigma,w_0)^{-1}$ for a positive $\nu$.
In view of Proposition 4.8, we now conclude

\medskip\noindent
{\bf Proposition 4.12}.
{\it The standard module $I(\nu,\sigma)=J(\nu,\sigma)$, i.e., $I(\nu,\sigma)$ is irreducible, if and only if $L(1,\tilde r\cdot\phi_\nu)^{-1}\neq 0$}.

The same is true when $k$ is non--archimedean and $\sigma$ is an irreducible generic tempered representation of $M(k)$ defined by a unitary character $\mu$ of $T(k)$ as in (3.22), (3.25) and (3.26), through induction.
This is a case of the unramified version of Proposition 4.12 for which LLC is automatic.
We therefore have:

\medskip\noindent
{\bf Proposition 4.13}.
{\it Assume $k$ is non--archimedean, $\sigma$ is generic, tempered and unramified whose Langlands parameter is $\phi$.
Then the standard module $I(\nu,\sigma)$ is irreducible if and only if $L(1,\tilde r\cdot\phi_\nu)^{-1}\neq 0$}.

Finally, we can state the following general result.
It is under the assumption of equality of $L$--functions through the local Langlands conjecture as discussed earlier.

\medskip\noindent
{\bf Proposition 4.14}.
{\it Assume $k$ is non--archimedean and $\sigma$ is an irreducible generic tempered representation of $M(k)$ whose Langlands parameter is $\phi$ for which
$$
L(s,r\cdot\phi_\nu)=L(s,\sigma_\nu,\tilde r),
$$
where $r$ is the adjoint action of ${}^L\!M$ on ${}^L\!\frak n$ and the $L$--function on the right is the one defined in \cite{48}.
Then the standard module $I(\nu,\sigma)$ is irreducible if and only if $L(1,\tilde r\cdot\phi_\nu)^{-1}\neq 0$}.

\begin{center}\bf 5. Proof of the Main Results\end{center}

In this section we prove Theorem 4.1 which we now generalize even further.

\medskip\noindent
{\bf Theorem 5.1}.
{\it Assume the validity of the local Langlands conjecture for every proper Levi subgroup $M$ of $G$ to the extent that every irreducible generic tempered representation $\sigma$ of $M(k)$ is parametrized by a homomorphism $\phi$ from $L_k$ to ${}^L\!M$ with a bounded image in $\widehat M$ such that
\[
L(s,r\cdot\phi_\nu)=L(s,\sigma_\nu,\tilde r),\tag{5.2}
\]
where $r$ and $\nu$ are as in Section 4.
Let $\psi\in\Psi (G/k)$ and let $\Pi(\phi_\psi)$ be the packet attached to $\phi_\psi$ defined by (3.6).
Suppose $\Pi(\phi_\psi)$ has a generic member.
Then $\phi_\psi$ is tempered}.

\medskip\noindent
{\bf Corollary 5.3}.
{\it Assume $k=\bR($or $\bC)$.
Let $\psi\in\Psi (G/\bR)$.
Then every generic member of $\Pi(\phi_\psi)$ is tempered.
In particular, if $\Pi(\phi_\psi)$ has a generic member then it is a tempered $L$--packet}.

\medskip\noindent
{\it Proof}.
This follows immediately from Theorem 5.1 since the full local Langlands conjecture for real groups is a theorem in \cite{37}.
The equality (5.2) is proved in \cite{46}.

Since unramified representations always satisfy LLC, the following corollary of Theorem 5.1 is immediate.

\medskip\noindent
{\bf Corollary 5.4}.
{\it Theorem 4.1 is valid}.

\medskip\noindent
{\it Proof of Theorem 5.1}.
Let $\pi$ be a generic member of $\Pi(\phi_\psi)$.
Assume $\pi=J(\nu,\tau)$, $\tau\in\Pi_\phi^M$ as in Proposition 3.31.
Then by Theorem 4.2, $I(\nu,\tau)$ is irreducible.

It now follows from Proposition 4.14 that
\[
L(1,\tilde r\cdot\phi_\psi)^{-1}\neq 0.\tag{5.5}
\]

Now assume the $\rho$--component of $\psi$ as defined in (3.8) is non--trivial, i.e., the $L$--packet $\Pi(\phi_\psi)$ is non--tempered.
Then by Proposition 3.31 and equation (3.20) and statement (3.21) the restriction $\tilde r_{\psi,\widehat\alpha}\colon\!\!=\tilde r\cdot \phi_\psi|X_{\widehat\alpha}$ acts like $w\mapsto |w|^{-1}$.
Consequently
\begin{align}
L(1,\tilde r\cdot\phi_\psi|X_{\widehat\alpha})^{-1}&=\zeta_k (0,|w|\cdot\tilde r_{\psi,\widehat\alpha}(w))^{-1}\nonumber\\
&= \zeta_k (0,\bold 1)^{-1}\tag{5.6}\\
&= 0,\nonumber
\end{align}
where $\zeta_k(s,\chi)$ denotes the Artin (or Hecke) $L$--function attached to a character $\chi$ of $L_k$.
This contradicts (5.5).
Thus $\rho=1$ and $\Pi(\phi_\psi)$ is tempered.
This completes the proof.

\medskip\noindent
(5.7)\ There are a number of cases where the local Langlands conjecture is proved \cite{19,22,24}.
The cases include certain cases of classical groups.
Consequently, our Theorem 5.1 is valid with no assumptions in those cases.
In particular, our proof of Theorem 5.1 gives a new proof of these results originally proved in \cite{9} and \cite{39}, without appealing to classification theorems.
The unramified case (Theorem 4.1) can also be proved using the results in \cite{40}.

\begin{center}\bf 6.\ Ramanujan Conjecture\end{center}

We now assume $k$ is a number field and $G$ is a quasisplit connected reductive algebraic group over $k$.
Let $\pi=\otimes_v\pi_v$ be a cuspidal automorphic representation of $G(\bA_k)$.
Going back to the discussions in Section 3, we now assume the following statement from a generalization of Arthur's $A$--packet conjecture due to Clozel (Conjecture 2A of \cite{m}).
In fact, our Conjecture 6.1 came out of our discussions with Arthur.
On the other hand, after consulting Clozel's article \cite{m} which was suggested to us by Dihua Jiang after a talk given on these results, it became clear that the question was already confronted by Clozel who then conjectured it as Conjecture 2A in \cite{m} by stating it as:\ {\it Suppose $\psi_v=(\phi_v,\rho_v)$ and $\phi_v|I'_{k_v}\equiv 1$.
Then the unramified members of $\Pi(\psi_v)$ are precisely those in $\Pi(\phi_{\psi_v})$.}

\medskip\noindent
{\bf Conjecture 6.1 (Arthur; Clozel, Conjecture 2A of \cite{m})}.
{\it For almost all finite primes $v,\pi_v\in\Pi (\phi_{\psi_v})$, where $\psi_v$ is the Arthur parameter of $\pi_v$}.

We now assume further that $\pi$ is locally generic.
Then each $\pi_v$ is generic with respect to a generic character $\chi_v$ of $U(k_v)$.
The characters $\chi_v$ may or may not be a local component of a global character of $U(k)\backslash U(\bA_k)$.
Appealing to Theorem 4.1 we conclude immediately that

\medskip\noindent
{\bf Theorem 6.2}.
{\it Assume Conjecture 6.1.
Let $\pi=\otimes_v\pi_v$ be a locally generic cuspidal automorphic representation of $G(\bA_k)$.
Then $\pi_v$ is tempered for almost all places $v$ of $k$}.

\medskip\noindent
{\it Remark 6.3}.
Since globally generic representations are locally generic, one may drop ``locally'' from our statement in this case.

Note that one can try to include more places in the tempered set whenever one can apply Theorem 5.1, but instead we will make the following conjecture for which we will produce some evidence.

\medskip\noindent
{\bf Conjecture 6.4}.
{\it Let $\pi=\otimes_v\pi_v$ be a cuspidal automorphic representation of $G(\bA_k)$.
Assume $\pi_v$ is tempered for almost all or even infinitely many places.
Then $\pi$ is tempered, i.e., $\pi_v$ is tempered for all $v$}.

\medskip\noindent
{\bf Corollary 6.5}.
{\it Assume Conjectures 6.1 and 6.4.
Let $\pi$ be a locally generic cuspidal automorphic representation of $G(\bA_k)$.
Then $\pi$ is tempered}.

Conjecture 6.4 is clearly valid for cuspidal representations attached to representations of the Galois group $\Gamma_k$, and using its validity for gr\"ossencharacters, for those parametrized by admissible homomorphisms from $W_k$ into $^LG$, whenever they exist.
One hopes that heuristic reasons can be given for the validity of Conjecture 6.5 in general, if one adopts the formalism of the global Langlands group $L_k$
(cf.~\cite{5,32}) and appeals to Arthur's global $A$--packet conjecture discussed in Section 3 (cf.~\cite{2,3,4}).
A much stronger conjecture is due to Clozel (Conjecture 4 of \cite{m}).

As another piece of evidence, we should mention the cases of classical groups.
Assuming the Ramanujan conjecture for $GL_N$, i.e., that all the unitary cuspidal representations of $GL_N(\bA_k)$ are tempered, one can show that the same is true for the generic spectrum of all the quasisplit classical groups.
This follows from the functorial transfer of generic cuspidal representations of such groups to appropriate $GL_N$ which was established in \cite{11,12,13,28}.
In fact, one knows that the transfers are always isobaric sums of unitary cuspidal representations of possibly smaller $GL_N$--groups and are therefore all tempered, assuming the Ramanujan conjecture for $GL_N$.
Thus, if $\pi=\otimes_v\pi_v$ is a globally generic cuspidal representation of $G(\bA_k)$, where $G$ is a quasisplit classical group over $k$, let $\Pi=\otimes_v\Pi_v$ be its functorial transfer to $GL_N(\bA_k)$.
We then conclude that $\pi_v$ is tempered for almost all $v$ if $\Pi$ satisfies the Ramanujan conjecture.
This extends to all $v$ as the work in Section 7 of \cite{12} shows; only a tempered $\pi_v$ can transfer to a tempered $\Pi_v$.
We refer to Section 10 of \cite{12} for details.
We collect this as 

(6.6).\ {\it Assume the validity of the Ramanujan conjecture for $GL_N(\bA_k)$.
Then globally generic cuspidal representations of $G(\bA_k)$ are all tempered}.

(6.7)\ Finally, we note that Corollary 6.5 in fact implies the Ramanujan conjecture for $GL_N(\bA_k)$ since cuspidal representations of $GL_N(\bA_k)$ are (locally) generic.

(6.8)\ We now address Conjecture 2.10 (and 6.1 and 6.4) in the case of classical groups by sketching a proof.
Here we start with assuming the Ramanujan--Selberg conjecture for $GL_N(\bA_k)$.
Arthur's upcoming book \cite{7} then allows us to define tempered packets for classical groups by transfering automorphic representations from $G(\bA_k)$ through functoriality to $GL_N(\bA_k)$.
In view of Corollary 6.5 all the generic representations are tempered and in particular so are globally generic ones.
On the other hand, by the automorphic descent of Ginzburg--Rallis--Soudry \cite{18,x}, globally generic representations of $G(\bA_k)$ are parametrized by certain self--dual isobaric automorphic representations of $GL_N(\bA_k)$ and the parametrization map is the inverse of the functoriality map from globally generic cusp forms to $GL_N(\bA_k)$ established in \cite{12}.
This transfer agrees with that of Arthur and in fact its image in $GL_N(\bA_k)$ is the image of the transfer of the whole tempered $L$--packet established in \cite{7}.
Since every tempered packet is transferred to a self--dual isobaric automorphic representation of $GL_N(\bA_k)$, all the generic ones in the same packet will transfer to the same representation on $GL_N(\bA_k)$.
In particular, the locally generic one in the packet must be the same as the globally generic one by uniqueness.

We finally remark that Conjectures 6.1 and 6.4 on which the validity of Corollary 6.5 relies, should be immediate consequences of the transfer established in \cite{11,12,13,28}.
In fact, a cuspidal representation $\pi$ whose local components are tempered at almost all places has no choice but to transfer to an isobaric automorphic representation of $GL_N(\bA_k)$ which will be tempered at every place by the Ramanujan conjecture for $GL_m(\bA_k)$, $1\leq m\leq N$.
This will then force all the components of $\pi$ to be tempered.
The validity of Conjecture 6.1 is readily available from the same transfer and the characterization of the residual spectrum for $GL_N(\bA_k)$ by Moeglin--Waldspurger \cite{41}.
We remark that this also gives a sketch of an argument for the validity of Conjecture 24.2 of \cite{16}.

\parskip=0pt
\bigskip\bigskip\noindent
Freydoon Shahidi

\noindent
Department of Mathematics

\noindent
Purdue University

\noindent
West Lafayette, IN\ 47907

\noindent
{\tt shahidi@math.purdue.edu}

\end{document}